\documentclass[reqno]{amsart}

\usepackage{amsmath}
\usepackage{amssymb}
\usepackage{amsfonts}
\usepackage{mathrsfs}
\usepackage{bbm}
\usepackage{algorithm2e}
\usepackage{MnSymbol}
\usepackage{caption}
\usepackage{extarrows}
\usepackage{cite}
\usepackage{graphicx}

\graphicspath{ {./Pictures/} }

\usepackage{todonotes}

\usepackage[tmargin=1.5in,bmargin=1.5in,lmargin=1.3in,rmargin=1.3in]{geometry}
\usepackage{microtype}

\newtheorem{thm}{Theorem}
\newtheorem{cor}{Corollary}
\newtheorem{op}{Open Problem}

\newtheorem{lem}{Lemma}
\theoremstyle{remark}

\title[Sharpness of the Phase Transition for the orthant model]{Sharpness of the Phase Transition for the orthant model}

\author{Thomas Beekenkamp}
\address{Mathematisches Institut, Ludwig-Maximilians-Universit\"at M\"unchen, Theresienstra\ss{}e 39, 80333 M\"unchen, Germany}
\email{Thomas.Beekenkamp@math.lmu.de}
\date{\today}

\begin{document}

\begin{abstract}
The orthant model is a directed percolation model on $\mathbb{Z}^d$, in which all clusters are infinite. We prove a sharp threshold result for this model: if $p$ is larger than the critical value above which the cluster of $0$ is contained in a cone, then the shift from $0$ that is required to contain the cluster of $0$ in that cone is exponentially small. As a consequence, above this critical threshold, a shape theorem holds for the cluster of $0$, as well as ballisiticity of the random walk on this cluster.
\end{abstract}

\maketitle

\section{Introduction and Main Result}
\begin{figure}
  \centering
    \includegraphics[width=0.75\textwidth]{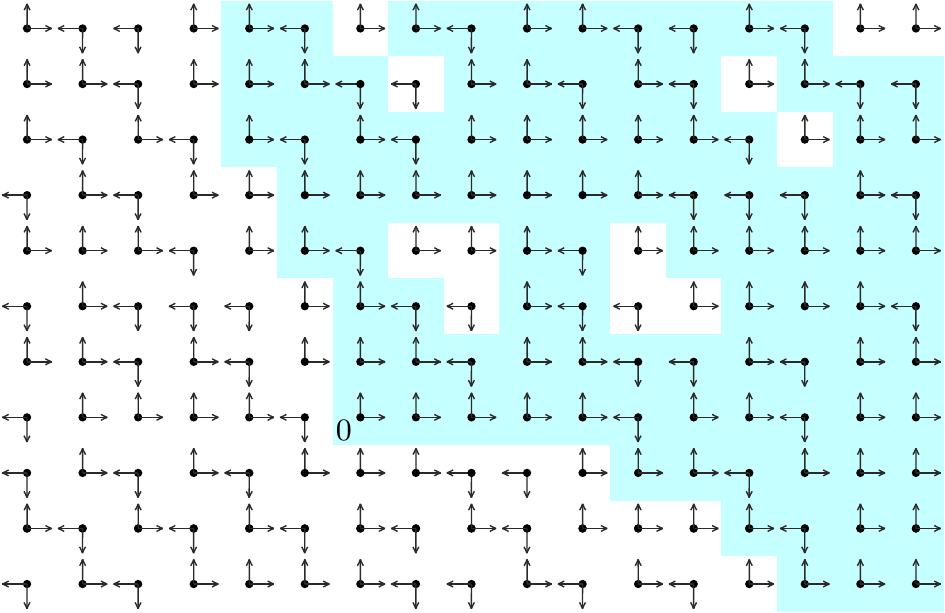}
    \caption{The orthant model on $\mathbb{Z}^2$. The cluster of the origin is shaded turquoise.}  
    \label{fig:orthant}
\end{figure}

We consider the orthant model on the directed graph $\mathbb{Z}^d$, $d\geq 2$, with nearest neighbour edges. This model is informally described as follows. Let $e_1,\dots,e_d$ be the standard unit basis vectors of $\mathbb{R}^d$. We set $\mathcal{E}_+:=\{e_1,\dots, e_d\}$, and $\mathcal{E}_-:=\{-e_1,\dots, -e_d\}$, as well as $\mathcal{E}=\mathcal{E}_+\cup\mathcal{E}_-$. A vertex $v\in \mathbb{Z}^d$ is connected to the vertices $v+e$ for all $e\in \mathcal{E}_+$ with a directed edge with probability $p$, independently of the other vertices. Otherwise, so with probability $1-p$, the vertex $v$ is connected to the vertices $v+e$ for all $e\in \mathcal{E}_-$. This model was introduced by Holmes and Salisbury \cite{Holmes2014a,Holmes2014}. The model is shown in Figure \ref{fig:orthant} for $d=2$.

The random directed graph obtained in this way has the property that every vertex is in an infinite cluster, since, for example, either the edge in the direction $e_1$ or in the direction $-e_2$ is always available. Therefore, there is no classical percolation phase transition in this model, where for small $p$ there are only finite clusters, and for large $p$ an infinite cluster exists. Instead, for small $p$ the clusters will have a tendency to move in the directions of $\mathcal{E}_-$, and for large $p$ a tendency in the directions of $\mathcal{E}_+$. In order to make this notion precise, we introduce a cone oriented in the direction $\mathbf{1}:=e_1+\dots+e_d$. For $0\leq \eta\leq 1$, we define the convex cone
\[
\mathcal{K}_{\eta}\:=\Big\{x\in\mathbb{R}^d\::\: x\cdot \mathbf{1}\geq \eta\|x\|_1\}.
\]

For $v,w\in\mathbb{Z}^d$, we say that $v\longrightarrow w$, whenever there is a directed path from $v$ to $w$. Note that this is not a symmetrical event, since we are working with a directed graph. Furthermore, for $A\subset \mathbb{Z}^d$, we say that $v\longrightarrow A$, whenever there exists $w\in A$ such that $v\longrightarrow w$. For $v\in \mathbb{Z}^d$, let $\mathcal{C}_v$ denote the forward cluster of $v$, i.e.,
\[
\mathcal{C}_v:=\{w\in \mathbb{Z}^d\::\: v\longrightarrow w\}.
\]
We can define the critical point above which $\mathcal{C}_0$ is contained in the cone with parameter $\eta$:
\[
\tilde{p}_c(\eta):=\inf\big\{p\::\: \mathbb{P}_p(\mathcal{C}_0\subset -n\mathbf{1}+\mathcal{K}_\eta \text{ for some }n\in \mathbb{N})=1\big\}.
\]
Note that this critical point is increasing in $\eta$, and that $\tilde{p}_c(1)=1$, so that the values $\eta\in[0,1)$ are of interest. In fact, Holmes and Salisbury \cite{Holmes2017} have proven that $p_c^{\text{OSP}}\leq\tilde{p}_c(\eta)<1$ for all $0<\eta<1$, where $p_c^{\text{OSP}}$ is the critical parameter for oriented site percolation on the triangular lattice. Furthermore, from considerations later in this section, it follows that $\tilde{p}_c(0)>0$. Therefore, $\tilde{p}_c(\eta)$ is non-trivial for $0\leq\eta<1$. In this paper, we will prove the following result.
\begin{thm}\label{thm}
Consider the orthant model on $\mathbb{Z}^d$ with parameter $p$. Let $0\leq \eta<1$ and suppose $p>\tilde{p}_c(\eta)$. Then, there exists a constant $c_p>0$, such that for all $n\in \mathbb{N}$,
\[
\mathbb{P}_p\big(0\longrightarrow (-n\mathbf{1}+\mathcal{K}_{\eta})^c\big)\leq \exp(-c_pn).
\]
\end{thm}
The above result for $d=2$ was proven by Holmes and Salisbury \cite{Holmes2014a} by making a connection with oriented site percolation on the triangular lattice. The result of Theorem \ref{thm} is known as a sharp threshold, or as a sharp phase transition. This type of result has been proven in a variety of models, most notably, and initially, for Bernoulli bond percolation in the 80's by Menshikov \cite{Menshikov1986} and Aizenman and Barsky \cite{AizenmanBarsky1987}. More recently, a revolutionary technique using the OSSS inequality was developed by Duminil-Copin, Raoufi and Tassion to prove sharp thresholds in models with more complexities \cite{Duminil-Copin2017,Duminil-Copin2018,Duminil-Copin2019,Duminil-Copin2020}.  These models include the random cluster model, Voronoi percolation, Boolean percolation, Gaussian fields, the corrupted compass model and the Widom-Rowlinson model \cite{Muirhead2020,Dereudre2018,Beekenkamp2020,Hutchcroft2020}.

As a consequence of the above sharp threshold result, we can prove a shape theorem for $\mathcal{C}_0$ above $\tilde{p}_c:=\lim_{\eta \downarrow 0}\tilde{p}_c(\eta)$. This critical point can also be written as
\[
\tilde{p}_c=\inf\big\{p\::\: \exists \eta>0 \text{ s.t. }\mathbb{P}_p(\mathcal{C}_0\subset -n\mathbf{1}+\mathcal{K}_\eta \text{ for some }n\in \mathbb{N})=1\big\}.
\]
A shape theorem for the orthant model was first proven by Holmes and Salisbury \cite{Holmes2019a} for large $p$. Using Theorem \ref{thm}, we can extend this result to all $p>\tilde{p}_c$. In order to state the shape theorem, we introduce for $u\in \mathbb{Z}^d$
\[
\beta_n(u):=\inf\{k\in \mathbb{Z}\::\: k\mathbf{1}+nu\in \mathcal{C}_0\}.
\]
Furthermore, let $\Lambda_r:=\{v\in \mathbb{Z}^d\::\: \|v\|_\infty\leq r\}$ be the closed ball around $0$ with radius $r$ with respect to the $L^\infty$-norm. Borrowing the notation from \cite{Holmes2019a}, the shape theorem for the orthant model can be stated as follows.

\begin{cor}[Shape theorem]\label{thm:shape}
Let $p>\tilde{p}_c$. The following hold for the orthant model on $\mathbb{Z}^d$ with parameter $p$. 
\begin{enumerate}
\item For $u\in \mathbb{Z}^d$, there is a deterministic $\gamma(u)\in \mathbb{R}$ such that $\frac{\beta_n(u)}{n}\to \gamma(u)$, as $n\to \infty$, $\mathbb{P}_p$-almost surely.
\item This limit satisfies $\gamma(u+w)\leq \gamma(u)+\gamma(w)$, $\gamma(ru)=r\gamma(u)$, $\gamma(u+r\mathbf{1})=\gamma(u)-r$, for $u,w\in \mathbb{Z}^d$, and $r\in \mathbb{N}$. Furthermore, $\gamma$ is symmetric under permutation of coordinates, $\gamma(u)\geq 0$ if $u\cdot\mathbf{1}\leq 0$, and $\gamma(u)\leq 0$ if $u$ lies in the positive orthant.
\item The limit $\gamma$ extends to a Lipschitz map $\mathbb{R}^d\to \mathbb{R}$ with these same properties, but for $r\in [0,\infty)$ and $u,w\in \mathbb{R}^d$.
\item The set $C:=\{z\in \mathbb{R}^d\::\:\gamma(z)\leq 0\}$ is a closed convex cone, which is symmetric under permutations of the coordinates, contains the positive orthant, and is contained in the half-space $\mathcal{K}_0=\{z\::\: z\cdot \mathbf{1}\geq 0\}$.
\item Let $\mathcal{C}^*_0:=\mathcal{C}_0+e_1 \mathbb{N}_0$, i.e., ``$\mathcal{C}_0$ with its holes filled in". It holds that $\frac{1}{n}\mathcal{C}^*_0\to C$, in the sense that for every $\varepsilon>0$ and every $r<\infty$, the following holds $\mathbb{P}_p$-a.s. for sufficiently large (random) n:
\[
\big(\Lambda_r\cap \frac{1}{n}\mathcal{C}^*_0\big)\,\subset\, \Lambda_\varepsilon+C,\quad\text{and}\quad (\Lambda_r\cap C) \,\subset\, \Lambda_r+\frac{1}{n}\mathcal{C}^*_0.
\]
\end{enumerate}
\end{cor}

To prove this theorem for all $p>\tilde{p}_c$, we modify the proof in \cite{Holmes2019a} by using Theorem \ref{thm} in the places where they require $p$ to be large. Another consequence of Theorem \ref{thm} is the ballisticity of the random walk on $\mathcal{C}_0$.
\begin{cor}[Ballisticity of the Random Walk]\label{thm:rw}
Consider the orthant model on $\mathbb{Z}^d$ with parameter $p>\tilde{p}_c$. Let $X_n$ be a simple random walk on $\mathcal{C}_0$ and let $P$ be the annealed law of this random walk (i.e., averaged over $\mathcal{C}_0$). Then $\tfrac{1}{n}X_n\to \mathbf{1}$ $P$-a.s. as $n\to \infty$, and
\[
\bigg(\frac{X_{\lfloor nt \rfloor}-\mathbf{1}nt}{\sqrt{n}}\bigg)_{t\geq 0} \Rightarrow (B_t)_{t\geq 0}, \quad \text{as }n\to \infty,
\]
weakly under $P$, where $(B_t)_{t\geq 0}$ is a $d$-dimensional Brownian motion with nonsingular covariance matrix $\Sigma$.
\end{cor}
This is Theorem 1.4 of \cite{Holmes2017} by Holmes and Salisbury applied to the orthant model. Their theorem is stated for more general models, and requires two conditions; one of which they show to hold for the orthant model with any value of $p$. The other condition is the existence of $\eta>0$ and $c>0$ such that $\mathbb{P}_p\big(0\longrightarrow (-n\mathbf{1}+\mathcal{K}_{\eta})^c\big)\leq \exp(-cn^\beta)$, for some $\beta>0$. By taking $\beta=1$ and assuming $p>\tilde{p}_c$, it follows from Theorem \ref{thm} that this condition holds for the orthant model with parameter $p$. Corollary \ref{thm:rw} is therefore an immediate consequence of combining Theorem \ref{thm} with Theorem 1.4 of \cite{Holmes2017}.

Despite the above results, the theoretical picture of the orthant model is still incomplete. We will use the remainder of this section to formulate two open questions for the model. The shape theorem and the ballisiticity of the random walk have now been shown to hold for $p>\tilde{p}_c:=\lim_{\eta \downarrow 0}\tilde{p}_c(\eta)$. A natural extension would be to prove these results for $p>\tilde{p}_c(0)$. This would follow from the continuity of $\tilde{p}_c(\eta)$.
\begin{op}
Consider the orthant model on $\mathbb{Z}^d$. The function $\eta\mapsto \tilde{p}_c(\eta)$ is continuous.
\end{op}
A critical value other than $\tilde{p}_c$ can be defined for the orthant model. In order to state this definition, we introduce for $v\in \mathbb{Z}^d$,
\[
L_v:=\inf\{k\in \mathbb{Z}\::\: v+ke_1\in \mathcal{C}_0\}.
\]
The critical value $p_c$ is defined as
\[
p_c:=\sup\{p\::\: L_0=-\infty \text{ a.s.}\}.
\] 
Holmes and Salisbury \cite{Holmes2019} have shown that this critical value is nontrivial, i.e., $0<p_c<1$. From the definitions of the critical values, it is clear that $p_c\leq \tilde{p}_c$. However, it is as of yet unclear that above $p_c$ there exists a cone with parameter $\eta>0$ that contains the forward cluster of $0$.
\begin{op}
Consider the orthant model on $\mathbb{Z}^d$. It holds that
\[
p_c=\tilde{p}_c.
\]
\end{op}
In order to prove this, perhaps it is most natural to first show that $p_c=\tilde{p}_c(0)$, and subsequently show the continuity of $\tilde{p}_c(\eta)$.

The sharp threshold result of Theorem \ref{thm} will be proven in Section \ref{sec:exp proof}, while some preliminaries required for this proof are introduced in Section \ref{sec:prelim}. The proof for Theorem \ref{thm:shape} is given in Section \ref{sec:shape proof}. 

\section{Preliminaries}\label{sec:prelim}
One difficulty in analysing the orthant model is the lack of monotonicity in $p$, i.e., a path from $v$ to $w$ might be lost if we increase $p$. To deal with this issue, we introduce the half-orthant model. In this model a vertex $v$ is always connected to $v+e$ for all $e\in \mathcal{E}_+$, whereas $v$ is connected to $v+e$, for all $e\in \mathcal{E}_-$, with probability $1-p$. This model is monotone in $p$, in the sense that $\mathbbm{1}\{v\longrightarrow w\}$ is monotonically decreasing in $p$. Let $\mathcal{C}^*_v$ denote the forward cluster of $v$ in the half-orthant model. The half-orthant model dominates the orthant model, in the sense that $\mathcal{C}_v\subseteq \mathcal{C}^*_v$, almost surely under a suitable coupling between the two models. For $v\in \mathbb{Z}^d$, we further define
\[
L^*_v:=\inf\{k\in \mathbb{Z}\::\: v+ke_1\in \mathcal{C}^*_0\}.
\]
From the domination it follows that $L_v\geq L^*_v$. However, it turns out that equality holds: $L_v^*=L_v$ for all $v\in \mathbb{Z}^d$ \cite[Thm. 1]{Holmes2019}. So, loosely speaking, if we only care about the leftmost boundary of $\mathcal{C}_0$, it does not matter if we consider the orthant model or the half-orthant model. This allows us to prove statements for the orthant model by making use of the monotonicty of the half orthant model. In light of this, we remark that the above definition of $\mathcal{C}_0^*$ coincides with the definition stated in Corollary \ref{thm:shape}. Furthermore, we note that $L_v<\infty$ for all $v\in \mathbb{Z}^d$, since $L^*_v<\infty$, but it might be the case that $L_v=-\infty$ for some $v\in \mathbb{Z}^d$. In fact, Holmes and Salisbury proved that if $L_v$ is finite for some $v\in \mathbb{Z}^d$, then it is finite for all $v\in \mathbb{Z}^d$. For $p<p_c$ it follows that $L^*_v=-\infty$ for all $v\in \mathbb{Z}^d$, and in this case, $\mathcal{C}^*_0=\mathbb{Z}^d$. On the other hand, if $p>p_c$, $L_v$ is finite for all $v\in \mathbb{Z}^d$  using the monotonicity of the half-orthant model. 

To prove Theorem \ref{thm}, it therefore suffices to work with the half-orthant model. We start by giving a formal description of this model. For $p\in[0,1]$, we consider the probability space $(\Omega, \mathcal{F}, \mathbb{P}_p)$, where
\[
\Omega=\{0,1\}^{\mathbb{Z}^d},
\]
the $\sigma$-algebra $\mathcal{F}$ is generated by the cylindrical events, and $\mathbb{P}_p$ is the product measure on $\Omega$ such that $\mathbb{P}_p(\omega_v=1)=p$ for all $v\in \mathbb{Z}^d$. From $\omega\in \Omega$ we obtain the edge configuration $\xi\subseteq\{(v,v+e)\::\:v\in \mathbb{Z}^d, e\in \mathcal{E}\}$ by adding the edge $(v,v+e)$ to the graph for all $e\in \mathcal{E}_+$, and for all $e\in \mathcal{E}_-$ whenever $\omega_v=0$.

For $v,w\in \mathbb{Z}^d$, we say $v\sim w$ when $v$ is a neighbour of $w$, i.e., whenever $w=v+e$ for some $e\in \mathcal{E}$. Furthermore, we say that $v\rightlsquigarrow w$, whenever $(v,w)\in \xi$. For $A\subset \mathbb{Z}^d$, we say that $v\xlongrightarrow{A}w$, whenever there is a path from $v$ to $w$ using only edges in $\xi$ with starting points in $A$. Note that $w$ does not have to be an element $A$ for this event to hold. For $A=\mathbb{Z}^d$ we use the shorthand notation $\{v\longrightarrow w \}:=\big\{v\xlongrightarrow{\mathbb{Z}^d}w\big\}$. Furthermore, the event $v\xlongrightarrow{A}v$ trivially holds for all $v\in \mathbb{Z}^d$, and all $A\subset \mathbb{Z}^d$.

The proof will make use of the OSSS inequality for Boolean functions $f:\Omega\to \{0,1\}$. In order to state this inequality, we introduce the influence of $v$ on $f$. This is defined as
\[
\text{Inf}_v:=\mathbb{P}_p(f(\omega)\neq f(\omega^{\oplus v})),
\]
where $\omega^{\oplus v}$ is given by
\[
(\omega^{\oplus v})_w=\begin{cases} 1-\omega_v & \text{if }w=v,\\
\omega_v & \text{if }w\neq v.
\end{cases}
\]
In other words, 
\[
\text{Inf}_v=\mathbb{P}_p(v\text{ is pivotal for the event }\{f=1\}).
\]

A decision tree $T$ is a random sequence of vertices $(v_0,v_1,\dots)$ that is build sequentially as follows. The tree starts by revealing the value of $\omega_{v_0}$, for the starting vertex $v_0$. Then, depending on the value of $\omega_{v_0}$, it chooses a vertex $v_1$ and reveals the value of $\omega_{v_1}$. This process continues until it has obtained enough information to determine the value of $f$, i.e., whenever the values of $\omega_v$ for unrevealed vertices $v$, cannot change $f$ any more. This leads to the definition of the revealment of $v$ by $T$:
\[
\text{Rev}_v(T):=\mathbb{P}_p(\omega_v \text{ is revealed by }T).
\]
The OSSS inequality states that for a Boolean function $f$ depending on finitely many variables and a decision tree $T$ that determines the value of $f$ we have
\[
\text{Var}(f)\leq \sum_{v\in \mathbb{Z}^d} \text{Inf}_v \text{Rev}_v(T).
\]
This inequality was proven by O'Donnell, Saks, Schramm and Servedio \cite{ODonnell2005}. A detailed exposition of Boolean functions has been written by O'Donnell \cite{ODonnell2014}, in which the proof of the OSSS inequality can also be found. Still, we will make use of Boolean functions $f$ that depend on infinitely many variables, so that we require an additional limit argument. The OSSS inequality can be generalised to Boolean functions with countable domains by a monotone convergence argument, provided  
\[
\mathbb{P}_p(f(\omega)\neq f(\omega_n))\to 0,\quad \text{ as }n\to \infty,
\]
where $\omega_n$ is given by 
\[
(\omega_n)_v=\begin{cases}\omega_v & \text{if } v \text{ is revealed by }T \text{ before step }n,\\
\tilde{\omega}_v & \text{otherwise},
\end{cases}
\]
where $\tilde{\omega}$ has law $\mathbb{P}_p$, and is independent of $\omega$. This is for example the case for decision trees that terminate in a finite number of steps on the set $\{f=1\}$. This generalisation of the OSSS inequality has been stated by Duminil-Copin, Raoufi and Tassion \cite{Duminil-Copin2019}.

\section{Proof of Theorem \ref{thm}}\label{sec:exp proof}
For $\eta\geq 0$ and $n\in \mathbb{N}$, we define the Boolean function 
\[
f_n:=\mathbbm{1}\{0\longrightarrow (-n\mathbf{1}+\mathcal{K}_{\eta})^c\}.
\]

\subsection{Exploration algorithm}
We now introduce decision trees that determine the value of $f_n$. A vital point in the proof is that we can uniformly bound the revealment of the vertices. If we only use one decision tree with a deterministic starting point, then the starting vertex will have revealment 1, so that we cannot find a nontrivial uniform bound on the revealment. Therefore, we will introduce the decision trees $T_k$, for  $1\leq k\leq n$, which all start at different vertices. In this way, we can average over $k$ and have a meaningful uniform bound on the revealment. The basic idea of the decision tree $T_k$ is that it explores the cluster of the boundary of $-k\mathbf{1}+\mathcal{K}_\eta$. If $0\longrightarrow (-n\mathbf{1}+\mathcal{K}_{\eta})^c $, this path must go through the boundary of the cone $-k\mathbf{1} +\mathcal{K}_\eta$, so that $T_k$ determines $f_n$. Furthermore, $T_k$ terminates in a finite number of steps when $f_n=1$.

\begin{figure}
  \centering
    \includegraphics[width=0.9\textwidth]{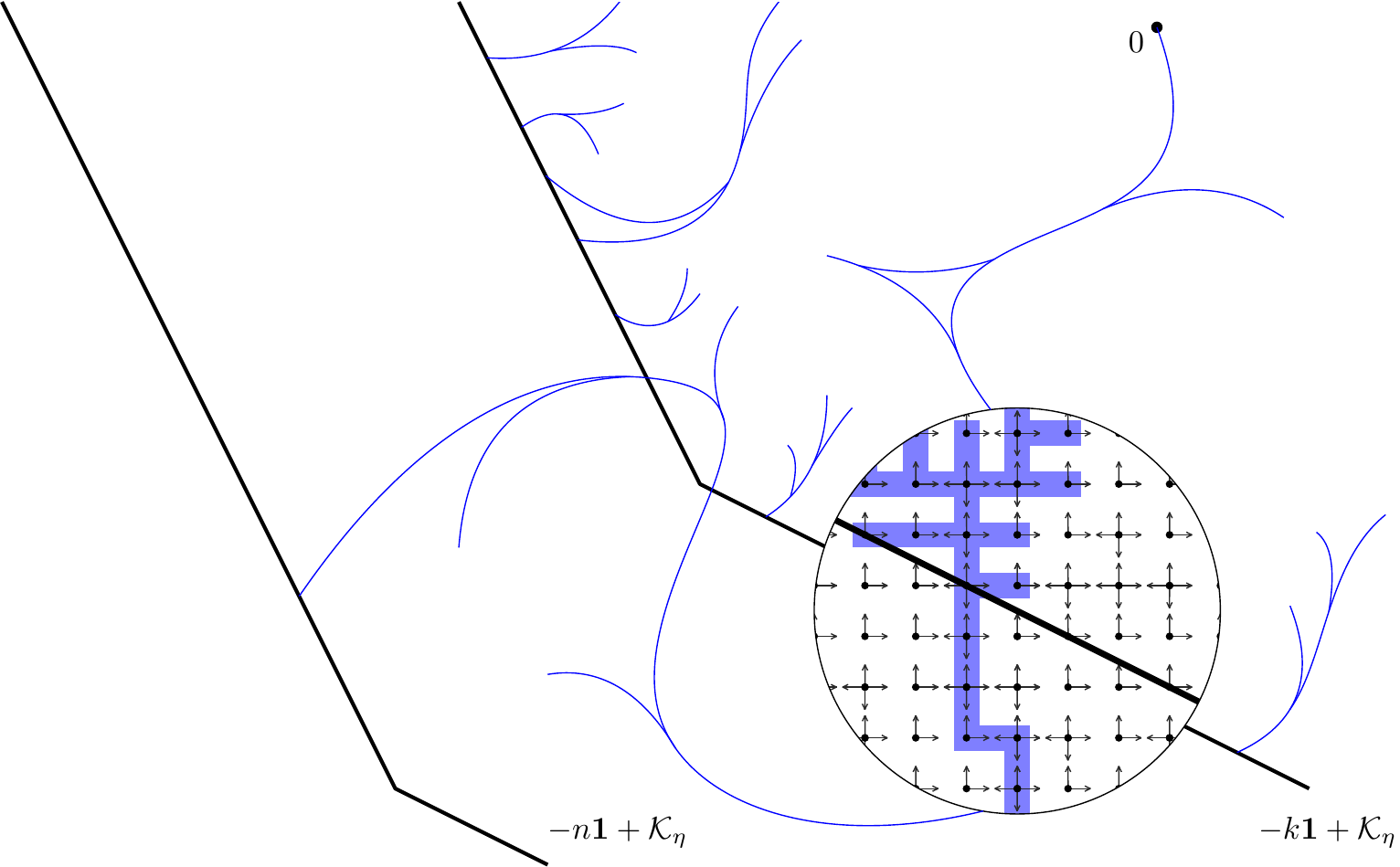}
    \caption{The algorithm $T_k$ exploring the cluster of $\partial(-k\mathbf{1}+\mathcal{K}_\eta)$ to find a path from $0$ to $(-n\mathbf{1}+\mathcal{K}_\eta)^c$. The blue vertices are revealed.}  
    \label{fig:algorithm}
\end{figure}

We will now describe the exploration algorithm of $T_k$ more precisely. We define the boundary and the outer boundary of the cone as
\begin{align*}
\partial (-k\mathbf{1} +\mathcal{K}_\eta)&:=\{v\in (-k\mathbf{1} +\mathcal{K}_\eta)\cap \mathbb{Z}^d \::\: \exists w\in (-k\mathbf{1} +\mathcal{K}_\eta)^c\cap\mathbb{Z}^d \text{ with }v\sim w\},\\
\partial^+ (-k\mathbf{1} +\mathcal{K}_\eta)&:=\{v\in (-k\mathbf{1} +\mathcal{K}_\eta)^c\cap \mathbb{Z}^d \::\: \exists w\in (-k\mathbf{1} +\mathcal{K}_\eta)\cap\mathbb{Z}^d \text{ with }v\sim w\}.
\end{align*}
The decision tree $T_k$ consists of two phases. In the first phase, $T_k$ explores the backward cluster of $\partial (-k\mathbf{1} +\mathcal{K}_\eta)$ inside the cone, that is, it explores the set $\{v\in -k\mathbf{1} +\mathcal{K}_\eta \::\: v\longrightarrow \partial (-k\mathbf{1} +\mathcal{K}_\eta)\}.$ When this is finished, the set of vertices 
\[
\left\{v\in \partial^+ (-k\mathbf{1} +\mathcal{K}_\eta)\::\:0\xlongrightarrow{-k\mathbf{1} +\mathcal{K}_\eta}v\right\}
\]
has been determined. In the second phase, the algorithm explores the forward clusters of these vertices. If for one of these vertices we find that $v\longrightarrow (-n\mathbf{1}+\mathcal{K}_\eta)^c$, then we also have $0\longrightarrow (-n\mathbf{1}+\mathcal{K}_\eta)^c$. A schematic visualisation of the algorithm is shown in Figure \ref{fig:algorithm}.

There is however one technical issue: since $f_n$ depends on the state of infinitely many vertices, it is possible that the algorithm gets stuck exploring inside $-k\mathbf{1}+\mathcal{K}_\eta$, and never gets to explore the forward clusters outside $-k\mathbf{1}+\mathcal{K}_\eta$. In order to deal with this, the decision tree operates in rounds, denoted by $i\in \mathbb{N}$. Recall that $\Lambda_r$ is the ball of radius $r$ around $0$ with respect to $L^\infty$-norm. In round $i$ we only explore inside $\Lambda_{i}$, so it is not possible to get stuck in any particular phase. Note that if $0\longrightarrow (-n\mathbf{1}+\mathcal{K}_\eta)^c$, there exists $i\in \mathbb{N}$ such that $0\xlongrightarrow{\Lambda_{i}} (-n\mathbf{1}+\mathcal{K}_\eta)^c$.

We denote by $\mathcal{R}$ the set of revealed vertices. Furthermore, we denote by $\mathcal{A}$ the set of active vertices for the first phase and by $\mathcal{B}$ the set of active vertices for the second phase. We start the algorithm by setting $\mathcal{A}:=\mathcal{A}_0:=\partial(-k\mathbf{1}+\mathcal{K}_\eta)$, and $\mathcal{B}:=\emptyset$. The pseudocode of $T_k$ is given in Algorithm \ref{alg}. We have to be careful when updating $\mathcal{A}$ in the first phase: note that by revealing $v$ it is possible that we create a new path $x\stackrel{\mathcal{R}}{\longrightarrow} \partial^+(-k\mathbf{1}+\mathcal{K}_\eta)$ for some $x\neq v$. Therefore, it is not sufficient to only consider $w\sim v$ for the update of $\mathcal{A}$. Instead, we add $w$ to $\mathcal{A}$ if and only if $w\not\in \mathcal{R}\cap \mathcal{B}$, and if there exists $x\in \mathcal{R}$ such that $x\stackrel{\mathcal{R}}{\longrightarrow} \partial^+(-k\mathbf{1}+\mathcal{K}_\eta)$.
\setlength{\algomargin}{0.9cm}
\begin{algorithm}
$i:=n$\;
$\mathcal{A}:= \partial (-k\mathbf{1} +\mathcal{K}_\eta)\cap\Lambda_n$\;
$\mathcal{B}:=\emptyset$\;
$\mathcal{R}:=\emptyset$\;
\While{$0 \stackrel{\mathcal{R}}{\not\longrightarrow} (-n\mathbf{1} +\mathcal{K}_\eta)^c$ }{
	\While{$\mathcal{A}\cap\Lambda_{i}\neq \emptyset$}{
	Take lexicographical minimal $v\in \mathcal{A}\cap\Lambda_{i}$\;
	Reveal $\omega_v$\;
	$\mathcal{R}:=\mathcal{R}\cup\{v\}$\;
	$\mathcal{A}:=\mathcal{A}\backslash \{v\}$\;
	$\mathcal{A}:=\mathcal{A}\cup \{w\in (-k\mathbf{1}+\mathcal{K}_\eta)\cap\mathbb{Z}^d\::\: w\not\in \mathcal{R}\cup\mathcal{B}, \exists x\in \mathcal{R}, x\sim w, \text{ s.t. }  x\stackrel{\mathcal{R}}{\longrightarrow} \partial^+(-k\mathbf{1}+\mathcal{K}_\eta) \}$\;
	$\mathcal{B}:=\mathcal{B}\cup \{w\in \mathbb{Z}^d\::\: w\not\in \mathcal{R},\exists x\in \partial^+(-k\mathbf{1}+\mathcal{K}_\eta)  \text{ s.t. }  0\stackrel{\mathcal{R}}{\longrightarrow} x, x\stackrel{\mathcal{R}}{\longrightarrow} w\}$\;
	$\mathcal{A}:=\mathcal{A}\backslash \{w\in \mathbb{Z}^d\::\: w\not\in \mathcal{R},\exists x\in \partial^+(-k\mathbf{1}+\mathcal{K}_\eta)  \text{ s.t. }  0\stackrel{\mathcal{R}}{\longrightarrow} x, x\stackrel{\mathcal{R}}{\longrightarrow} w\}$\;
	\lIf{$0 \stackrel{\mathcal{R}}{\longrightarrow} (-n\mathbf{1} +\mathcal{K}_\eta)^c$}{\Return{1}}
	}
	\While{$\mathcal{B}\cap\Lambda_{i}\neq \emptyset$}{
	Take lexicographical minimal $v\in \mathcal{B}\cap\Lambda_{i}$\;
	Reveal $\omega_v$\;
	$\mathcal{R}:=\mathcal{R}\cup\{v\}$\;
	$\mathcal{B}:=\mathcal{B}\backslash\{v\}$\;
	$\mathcal{B}:=\mathcal{B}\cup \{w\in \mathbb{Z}^d\::\: w\not\in \mathcal{R},\exists x\in \partial^+(-k\mathbf{1}+\mathcal{K}_\eta)  \text{ s.t. }  0\stackrel{\mathcal{R}}{\longrightarrow} x, x\stackrel{\mathcal{R}}{\longrightarrow} w\}$\;
	$\mathcal{A}:=\mathcal{A}\backslash \{w\in \mathbb{Z}^d\::\: w\not\in \mathcal{R},\exists x\in \partial^+(-k\mathbf{1}+\mathcal{K}_\eta)  \text{ s.t. }  0\stackrel{\mathcal{R}}{\longrightarrow} x, x\stackrel{\mathcal{R}}{\longrightarrow} w\}$\;
	\lIf{$0 \stackrel{\mathcal{R}}{\longrightarrow} (-n\mathbf{1} +\mathcal{K}_\eta)^c$}{\Return{1}}	
	}
	$i:=i+1$\;
}
\caption{The exploration algorithm $T_k$.}\label{alg}
\end{algorithm}
At the start of any iteration of the inner loops of the algorithm, the following hold for the active sets $\mathcal{A}$ and $\mathcal{B}$:
\begin{align*}
\mathcal{A}&\subseteq\mathcal{A}_0\backslash\mathcal{R}\cup\{v\in (-k\mathbf{1}+\mathcal{K}_\eta)\cap\mathbb{Z}^d \::\: v\not\in \mathcal{R}, \exists w\in \mathcal{R}\cap (-k\mathbf{1}+\mathcal{K}_\eta), w\sim v, \text{ s.t. } w\stackrel{\mathcal{R}}{\longrightarrow} \partial^+(-k\mathbf{1}+\mathcal{K}_\eta)\},\\
\mathcal{B}&=\{v\in (-n\mathbf{1}+\mathcal{K}_\eta)\cap \mathbb{Z}^d \::\: v\not\in \mathcal{R}, \exists w\in \partial^+(-k\mathbf{1}+\mathcal{K}_\eta)  \text{ s.t. }  0\stackrel{\mathcal{R}}{\longrightarrow} w, w\stackrel{\mathcal{R}}{\longrightarrow} v\}.
\end{align*}

\subsection{Bound on the revealment}
Let $\theta_n(p):=\mathbb{P}_p(f_n=1)$. Summing the OSSS inequality over $k$ gives
\begin{align}\label{eq:osss1}
n\theta_n(p)(1-\theta_n(p))\leq \sum_{v\in \mathbb{Z}^d}\text{Inf}_v \sum_{k=1}^n\text{Rev}_v(T_k).
\end{align}

We will now bound $\sum_{k=1}^n\text{Rev}_v(T_k)$ uniformly in $v$. Let $k_v$ be such that $v\in \partial(-k_v\mathbf{1}+\mathcal{K}_\eta)$. Suppose first that $k> k_v$. If $v$ is revealed by $T_k$ in the second phase, we have $0\stackrel{\mathcal{R}}{\longrightarrow} v$, so that in particular $0\longrightarrow \partial (-k\mathbf{1}+\mathcal{K}_\eta)$. On the other hand, if $v$ is revealed by $T_k$ in the first phase, there exists $w\sim v$ such that $w\longrightarrow \partial^+(-k\mathbf{1}+\mathcal{K}_\eta)$. Applying the union bound gives
\[
\sum_{k=1}^n \mathbbm{1}\{k_v<k\}\text{Rev}_v(T_k)\leq \sum_{k=1}^n \mathbbm{1}\{k_v<k\}\Big(\theta_k(p)+\sum_{w\sim v}\mathbb{P}_p(w\longrightarrow \partial^+(-k\mathbf{1}+\mathcal{K}_\eta))\Big).
\]
Let $d(w,-k\mathbf{1}+\mathcal{K}_\eta)$ be the distance between $w$ and $-k\mathbf{1}+\mathcal{K}_\eta$ in the $L^2$-norm. Then
\[
w-d(w,-k\mathbf{1}+\mathcal{K}_\eta)\mathbf{1}+\mathcal{K}_\eta\subseteq -k\mathbf{1}+\mathcal{K}_\eta,
\]
since $-k\mathbf{1}+\mathcal{K}_\eta$ is a convex cone, and $w-d(w,-k\mathbf{1}+\mathcal{K}_\eta)\mathbf{1}\in -k\mathbf{1}+\mathcal{K}_\eta$. Therefore, using translation invariance, it follows that
\begin{align}
\sum_{k=1}^n \mathbbm{1}\{k_v<k\}\text{Rev}_v(T_k)&\leq \sum_{k=1}^n \mathbbm{1}\{k_v<k\}\Big(\theta_k(p)+\sum_{w\sim v}\mathbb{P}_p(w\longrightarrow \partial(w -d(w, -k\mathbf{1}+\mathcal{K}_\eta)\mathbf{1}+\mathcal{K}_\eta))\Big)\nonumber\\
&=\sum_{k=1}^n\mathbbm{1}\{k_v<k\}\Big(\theta_k(p)+\sum_{w\sim v}\mathbb{P}_p(0\longrightarrow \partial( -d(w, -k\mathbf{1}+\mathcal{K}_\eta)\mathbf{1}+\mathcal{K}_\eta))\Big)\nonumber\\
&\leq\sum_{k=1}^n \mathbbm{1}\{k_v<k\}\theta_k(p)+2d\sum_{k=0}^{n}\mathbbm{1}\{k_v<k\}\mathbb{P}_p(0\longrightarrow \partial(-k\mathbf{1}+\mathcal{K}_\eta))\nonumber\\
&=\sum_{k=1}^n\mathbbm{1}\{k_v<k\}\theta_k(p)+2d\sum_{k=0}^{n}\mathbbm{1}\{k_v<k\}\theta_k(p).\label{eq: rev1}
\end{align}
Now suppose $k< k_v$, so $v\not\in -k\mathbf{1}+\mathcal{K}_\eta$. If $v$ is revealed, it holds that $0\stackrel{\mathcal{R}}{\longrightarrow} v$. In particular, we have $0\longrightarrow \partial^+(-k\mathbf{1}+\mathcal{K}_\eta)$. We find
\begin{align}
\sum_{k=1}^{n}\mathbbm{1}\{k_v>k\} \text{Rev}_v(T_k)\leq \sum_{k=1}^{n} \mathbbm{1}\{k_v<k\}\theta_{k}(p).\label{eq:rev2}
\end{align}
Combining (\ref{eq: rev1}) and (\ref{eq:rev2}) gives
\[
\sum_{k=1}^n \text{Rev}_v(T_k)\leq 1+\sum_{k=1}^{n}\theta_{k}(p)+2d\sum_{k=0}^{n}\theta_k(p)= (2d+1)\sum_{k=0}^{n}\theta_k(p).
\]
Writing $S_n:=\sum_{k=0}^n\theta_k(p)$, gives
\begin{align}\label{eq:osss2}
\sum_{v\in \mathbb{Z}^d}\text{Inf}_v  \geq \frac{1}{2d+1}\frac{n}{S_{n}}\theta_n(p)(1-\theta_n(p)).
\end{align}

\subsection{Analysis of the differential inequality}
We are now able to complete the proof of Theorem \ref{thm}. We can obtain a differential inequality by using Russo's formula. However, since $f_n$ depends on infinitely many vertices, $\theta_n(p)$ is not necessarily differentiable in $p$. Instead we have to work with the upper-right Dini derivative:
\[
D^+ \theta_n(p):= \limsup_{h\downarrow 0} \frac{\theta_n(p+h)-\theta_n(p)}{h}.
\]
Using the fact that $0\longrightarrow (-n\mathbf{1}+\mathcal{K}_{\eta})^c$ is a decreasing event, i.e., $f_n$ is a decreasing function of $\omega$, Russo's formula gives
\[
-D^+ \theta_n(p)\geq\sum_{v\in \mathbb{Z}^d}\mathbb{P}_p(v\text{ is pivotal for }\{f_n=1\})=\sum_{v\in \mathbb{Z}^d} \text{Inf}_v.
\]
This version of Russo's formula can be found in the book on Percolation by Grimmett \cite{Grimmett1999}. This is the point in the proof where use the monotonicity of the half-orthant model. Combining the above inequality with (\ref{eq:osss2}) gives
\begin{align}\label{eq:osss3}
-D^+ \theta_n(p)\geq\frac{1}{4d}\frac{n}{S_{n}}\theta_n(p)(1-\theta_n(p)),
\end{align}
where we use $2d+1\leq4d$ for simplicity. The rest of the proof consists of analysing the above diferential inequality. This analysis follows the line of Duminil-Copin, Raoufi and Tassion, but since it differs on several points, we choose to include it. We have to work with Dini derivatives instead of regular derivatives, and, more importantly, in our case we cannot give a simple lower bound on $1-\theta_n(p)$.

To analyse the differential inequality, we introduce the auxillary critical point
\[
\hat{p}_c(\eta):=\sup\left\{p\::\: \limsup_{n\to \infty} \frac{\log S_n(p)}{\log n}= 1\right\}.
\]
We will first show that $\hat{p}_c(\eta)\leq \tilde{p}_c(\eta)$, for $\eta\geq 0$. To prove this, we assume the contrary, and let $p\in(\tilde{p}_c(\eta), \hat{p}_c(\eta))$. Since $p>\tilde{p}_c(\eta)$, we can fix $N\in \mathbb{N}$, such that for all $n>N$ it holds that $\theta_{\lceil \sqrt{n}\rceil}(p)\leq 1/2$. We define $T_n(p):=\frac{2}{\log n}\sum_{k=\sqrt{n}}^n \frac{\theta_k(p)}{k}$. Taking the upper-right Dini derivative and using (\ref{eq:osss3}) gives
\begin{align}\label{eq:osss4}
-D^+ T_n &\geq \frac{1}{2d}\frac{1}{\log n} \sum_{k=\sqrt{n}}^n \frac{\theta_k(p)}{S_k}(1-\theta_k(p))\geq\frac{1}{4d}\frac{1}{\log n} \sum_{k=\sqrt{n}}^n \frac{\theta_k(p)}{S_{k}}\geq\frac{1}{4d}\frac{\log S_{n+1}-\log S_{\sqrt{n}}}{\log n},
\end{align}
where in the last inequality we used
\[
\frac{\theta_k(p)}{S_k}\geq \int_{S_{k}}^{S_{k+1}} \frac{1}{x}\,\mathrm{d}x=\log S_{k+1}-\log S_{k}.
\]
Now let $p_1\in(p,\hat{p}_c(\eta))$. We will integrate the differential inequality between $p$ and $p_1$ and use the following result regarding Dini derivatives: the Dini derivative of a decreasing function $f:[a,b]\to \mathbb{R}$ satisfies
\begin{align}\label{eq:dini}
f(b)-f(a)\leq \int_a^b D^+ f(x)\,\mathrm{d}x.
\end{align}
Applying this to $T_n(p)$ and using (\ref{eq:osss4}) gives
\[
T_n(p_1)-T_n(p)\leq \int_{p}^{p_1} D^+ T_n(s)\,\mathrm{d}s\leq -(p_1-p)\frac{1}{4d}\frac{\log S_{n+1}(p_1)-\log S_{\sqrt{n}}(p)}{\log n}.
\]
Furthermore, $T_n(p)$ converges to $\theta(p):=\lim_{n\to \infty}\theta_n(p)$ for $n\to \infty$, since
\[
\theta_{n}(p)=2\theta_n(p)\frac{\log n-\log \sqrt{n}}{\log n}\leq T_n(p) \leq 2\theta_{\lceil\sqrt{n}\rceil}(p)\frac{\log n-\log \sqrt{n}}{\log n}=\theta_{\lceil\sqrt{n}\rceil}(p).
\]
It follows that
\[
\theta(p_1)-\theta(p)\leq -(p_1-p)\frac{1}{4d}\limsup_{n\to \infty} \frac{\log S_{n+1}(p_1)-\log S_{\sqrt{n}}(p)}{\log n}.
\]
Since $p<\hat{p}_c(\eta)$, we have that $\limsup_{n\to \infty} \frac{\log S_n(p)}{\log n}= 1$ and the same holds for $p_1$, it follows that
\[
\limsup_{n\to \infty} \frac{\log S_{n+1}(p)-\log S_{\sqrt{n}}(p_1)}{\log n}=\limsup_{n\to \infty} \frac{\log S_{n+1}(p)}{\log n}-\frac{\log S_{\sqrt{n}}(p_1)}{2\log \sqrt{n}}= \frac{1}{2}.
\]
We conclude
\[
\theta(p)\geq \theta(p)-\theta(p_1)\geq\frac{p_1-p}{8d}>0,
\]
which contradicts $p>\tilde{p}_c(\eta)$, so that we have established that $\hat{p}_c(\eta)\leq \tilde{p}_c(\eta)$.

Now suppose $p>\tilde{p}_c(\eta)$, so that also $p>\hat{p}_c(\eta)$. Then there exists $\beta<1$ such that $S_n(p)\leq n^\beta$ and there exists $N\in \mathbb{N}$ such that $\theta_n(p)\leq \tfrac{1}{2}$ for all $n\geq N$. Combining this with (\ref{eq:osss3}) and using the chain rule for Dini derivatives gives
\[
D^+ \log \theta_n(p) \leq - \frac{1}{4d}n^{1-\beta}(1-\theta_n(p))\leq -\frac{1}{8d}n^{1-\beta},
\]
for $n>N$. Let $p_1:=(\tilde{p}_c(\eta)+p)/2$. Integrating the above inequality between $p_1$ and $p$ and using (\ref{eq:dini}) gives
\[
\log \theta_n(p)\leq \log \theta_n(p)-\log\theta_n(p_1) \leq -\frac{1}{8d}(p-p_1)n^{1-\beta}.
\]
It follows that
\[
\theta_n(p)\leq \exp\left(-\frac{1}{16d}(p-\tilde{p}_c(\eta))n^{1-\beta}\right).
\]
It remains to improve the above stretched exponential decay to proper exponential decay. From the stretched exponential decay it follows that $S(p):=\lim_{n\to \infty}S_n(p)<\infty$. Combining this fact with (\ref{eq:osss3}), and using that $\theta_n(p)\leq \tfrac{1}{2}$ for $n>N$, since $p>\hat{p}_c(\eta)$, gives
\[
D^+ \log \theta_n(p) \leq -\frac{1}{8dS}n.
\]
From here the proof is similar as for the stretched exponential decay, and we conclude
\[
\theta_n(p)\leq \exp\left(-\frac{1}{16dS(p)}(p-\tilde{p}_c(\eta))n\right).
\]
Theorem \ref{thm} now holds with 
\[
c_p:=\frac{1}{16dS(p)}(p-\tilde{p}_c(\eta))\wedge \sup\big\{C>0\::\: \theta_n(p)\leq \exp(-Cn) \text{ for all }n\leq N\big\}>0.
\]
\qed

\section{Proof of the Shape Theorem}\label{sec:shape proof}
To prove Corollary \ref{thm:shape}, we modify the proof of Holmes and Salisbury \cite{Holmes2019a} in the places where they require $p$ to be large. Their proof is structured in seven lemmas, two of which require a large $p$. The first of these is Lemma 1 of \cite{Holmes2019a}. This lemma asserts the existence of $\theta>1$, such that for every $\eta\in[0,1)$, there exists $p_0=p_0(\eta,d)<1$, such that for $p>p_0$, there exists $c_1>0$ such that $\mathbb{P}_p(0\longrightarrow (-n\mathbf{1}+\mathcal{K}_\eta)^c)\leq c_1 \theta^{-nd}$, for all $n\in \mathbb{N}$. In the remainder of their proof, this lemma is only used for the case $\eta=0$. Therefore, we can replace this lemma by Theorem \ref{thm}, and require $p>\tilde{p}_c$, instead of $p>p_0$.

The second lemma in the proof of Holmes and Salisbury which require large $p$ is Lemma 5 of \cite{Holmes2019a}. We will prove this lemma for $p>\tilde{p}_c$, instead of for large $p$, using Theorem \ref{thm}. To state this lemma, we let $u\in \mathbb{Z}^d\backslash \mathbb{Z}\mathbf{1}$, and fix $v\in \mathbb{R}^d$ such that $u\cdot v>0$ and $v\cdot \mathbf{1}=0$. We define the slab
\[
\Lambda_{u,v}(m,n):=\{z\in \mathbb{Z}^d\::\: mu\cdot v\leq z\cdot v< nu\cdot v\}.
\]
We are interested in the following three events. Let $A_n'(M)$ be the event there exists a path starting in 0 and ending in a point $k\mathbf{1}+nu$ with $k<n\gamma(u)$ that hits $\Lambda_{u,v}(-\infty, -M)$. Let $A_n''(M)$ be the event there exists a path starting in 0 and ending in a point $k\mathbf{1}+nu$ with $k<n\gamma(u)$ that hits $\Lambda_{u,v}(M+n, \infty)$. Lastly, let $\hat{A}_n$ be the event that there is a path starting in 0 and ending in some point $k\mathbf{1}$, with $k<0$, and reaches $\Lambda_{u,v}(n,\infty)$. We will prove the following lemma regarding these events:

\begin{lem}
Let $p>\tilde{p}_c$. There exists $c>0$, such that $\mathbb{P}_p(A_n'(\lfloor cn \rfloor)\text{ i.o.})=\mathbb{P}_p(A_n''(\lfloor cn \rfloor)\text{ i.o.})=\mathbb{P}_p(\hat{A}_n(\lfloor cn \rfloor)\text{ i.o.})=0$.
\end{lem}

\begin{figure}
  \centering
    \includegraphics[width=0.95\textwidth]{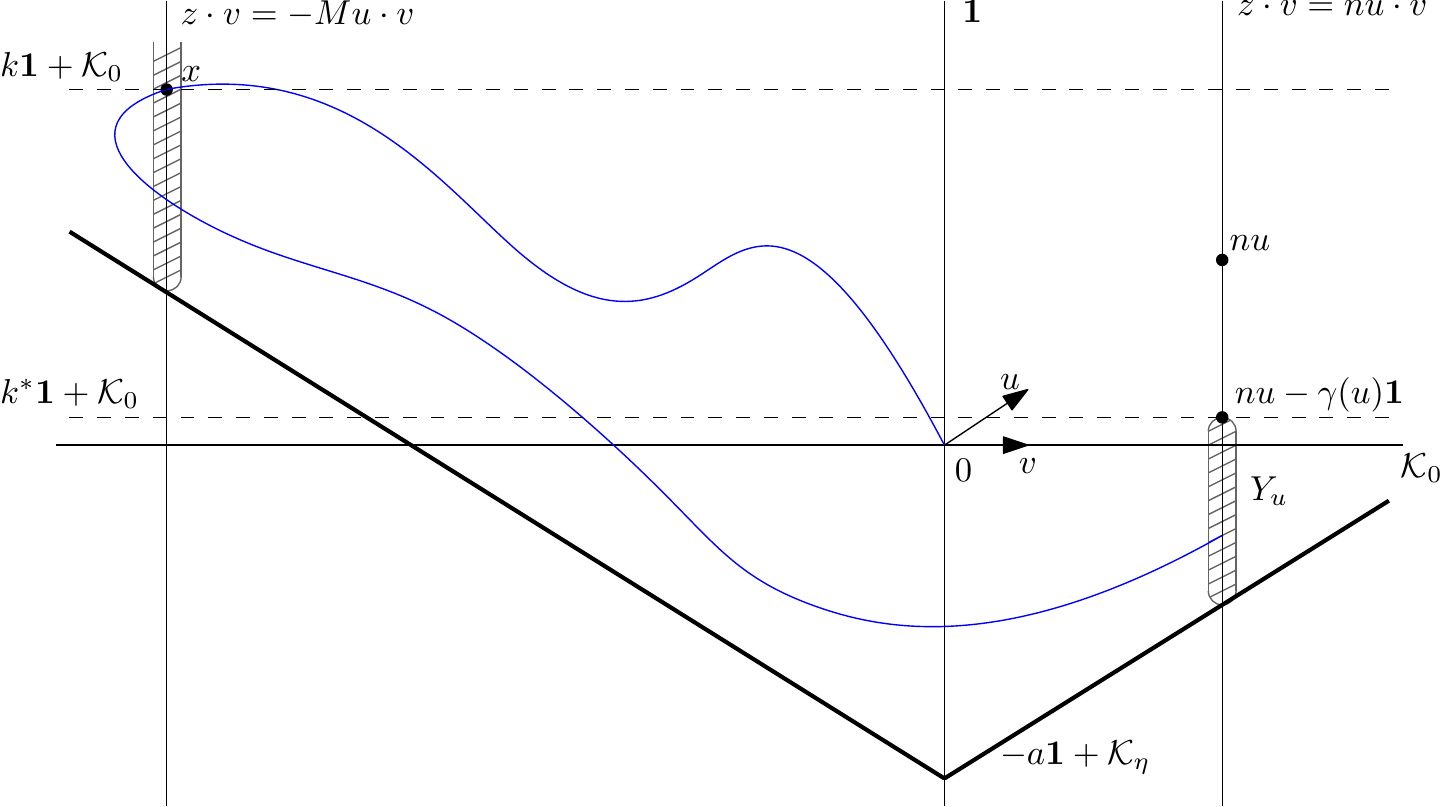}
    \caption{When the event $A_n'(M)$ occurs, there exists a path from $0$ to $Y_u$ going through the shaded region}  
    \label{fig:lem5}
\end{figure}

We will prove the above lemma for the event $A_n'(\lfloor cn \rfloor)$, the other two events can be proven similarly. The event $A_n'(\lfloor cn \rfloor)$ is shown in Figure \ref{fig:lem5}. Let $p>\tilde{p}_c$. By the definition of this critical point there exists $\eta>0$ such that $p>\tilde{p}_c(\eta)$. We fix such an $\eta$. Let $c>0$, and let $M:=M(n):=\lfloor cn\rfloor$. We will choose the precise value of $c$ later on. Let $a>0$ and suppose $\mathcal{C}_0^*\subseteq -a\mathbf{1}+\mathcal{K}_\eta$. If $A_n'(M)$ occurs, there exists $x\in \mathbb{Z}^d$ satisfying
\begin{equation}\label{eq:lem5sytem}
\begin{cases}
(x+a\mathbf{1})\cdot \mathbf{1}\geq \eta\|x+a\mathbf{1}\|_1,\\
\quad\quad \quad x\cdot v =-M u\cdot v,
\end{cases}
\Longrightarrow
\begin{cases}
x\cdot \mathbf{1}\geq -2da +\eta\|x\|_1,\\
 x\cdot v =-M u\cdot v,
\end{cases}
\end{equation}
such that $0\longrightarrow x$, and $x\longrightarrow y$, with $y=k\mathbf{1}+nu$ for some $k<n\gamma(u).$ Since the $L^1$-norm is equivalent to the $L^2$-norm, and since the $L^2$-norm is invariant under an orthonormal basis change, it follows from the above equation that $\|x\|_1\geq c_0 M $, for some constant $c_0=c_0(u,v)>0$. Combining this with the above inequality gives
\[
x\cdot\mathbf{1}\geq -2da+c_0\eta M.
\]
We define the set
\[
Y_u:=\{y\in \mathbb{Z}^d\::\: y=k\mathbf{1}+nu,\text{ with } k<n\gamma(u)\}.
\]
We use Theorem \ref{thm} and the union bound to obtain
\[
\mathbb{P}_p(A_n'(M))\leq \exp(-c_p a)+\!\sum_{k=-2da+c_0\eta M}^\infty\mathbb{P}_p(\exists x\in \Lambda_{u,v}(-\infty,-M)\::\:x\cdot\mathbf{1}=k,\, 0\longrightarrow x,\, x\longrightarrow Y_u).
\]
We define $k^*:=k^*(n):=n(d\gamma(u)+u\cdot\mathbf{1})$. With this choice, it follows that $y\cdot\mathbf{1}<k^*$ for all $y\in Y_u$, and all $n\in \mathbb{N}$. Suppose $x\cdot\mathbf{1}=k$, and $x\longrightarrow Y_u$, then it follows, that $x-k\mathbf{1}\longrightarrow (-(k-k^*)\mathbf{1}+\mathcal{K}_{0})^c$. We now fix
\[
c:=\Big(\frac{d\gamma(u)+u\cdot \mathbf{1}}{c_0\eta}+1\Big)\vee 1.
\]
Then, for $k=-2da+c_0\eta M+k'$, with $k'\geq 0$ it holds, that
\[
k-k^*= -2da +c_0\eta \Big\lfloor n\Big(\frac{d\gamma(u)+u\cdot \mathbf{1}}{c_0\eta}+1\Big)\Big\rfloor+k'-n(d\gamma(u)+u\cdot\mathbf{1})\geq -2da+c_0\eta(n-1)+k'=:f(n,k').
\]
It follows, that
\begin{align*}
\mathbb{P}_p(A_n'(M))\leq \exp(-c_p a)+\sum_{k'=0}^\infty\mathbb{P}_p\big(\exists x\in \Lambda_{u,v}(-\infty,-M)\:&:\:x\cdot\mathbf{1}=-2da+c_0\eta M+k',\\
&\quad x-x\cdot\mathbf{1}\longrightarrow -f(n,k')\mathbf{1}+\mathcal{K}_{0})^c\big).
\end{align*}
Combining $x\cdot\mathbf{1}=-2da+c_0\eta M+k'$ with (\ref{eq:lem5sytem}), shows that
\[
\|x\|_1\leq c_0 M+\frac{k'}{\eta}.
\]
Using another union bound, translation invariance, and Theorem \ref{thm}, we find
\begin{align}\label{eq:lem5end}
\mathbb{P}_p(A_n'(M))&\leq \exp(-c_p a)+\!\sum_{k'=0}^\infty \Big|\Big\{x\in \mathbb{Z}^d\::\: \|x\|_1\leq c_0 M+\frac{k'}{\eta}\Big\}\Big|\mathbb{P}_p(0\longrightarrow (-f(n,k')\mathbf{1}+\mathcal{K}_{0})^c)\nonumber\\
&\leq \exp(-c_p a)+\!\sum_{k'=0}^\infty \Big(2c_0 cn+2\frac{k'}{\eta}\Big)^d\exp(-c_p f(n,k'))
\end{align}
We now take 
\[
a:=a(n):=\frac{c_0\eta}{4d}n,
\]
so that
\[
f(n,k')=-2da+c_0\eta(n-1)+k'= \tfrac{1}{2}n-c_0\eta+k'.
\]
A careful examination of (\ref{eq:lem5end}) shows that the sum over $k'$ converges for all $n\in \mathbb{N}$, and that the result is summable with respect to $n$, so that by the Borel-Cantelli lemma $\mathbb{P}_p(A_n'(\lfloor cn \rfloor)\text{ i.o.})=0$. The same result can be proven similarly for the events $A_n''(\lfloor cn \rfloor)$ and $\hat{A}_n$, and we omit the proof.
\qed

\section*{Acknowledgments}
I thank Matija Pasch for insightful discussions on the topic, as well as for useful comments on the manuscript.
\bibliographystyle{abbrv}
\bibliography{om}

\end{document}